\documentclass[reqno,10pt,centertags,draft]{amsart}
\usepackage{amssymb}
\newcommand{\bbN}{{\mathbb{N}}}
\newcommand{\bbR}{{\mathbb{R}}}

\newcommand{\bbC}{{\mathbb{C}}}

\newcommand{\calB}{{\mathcal B}}
\newcommand{\calC}{{\mathcal C}}
\newcommand{\calD}{{\mathcal D}}

\newcommand{\calH}{{\mathcal H}}

\newcommand{\calK}{{\mathcal K}}

\newcommand{\calW}{{\mathcal W}}

\newcommand{\no}{\nonumber}
\newcommand{\lb}{\label}

\newcommand{\bi}{\bibitem}

\newcommand{\beq}{\begin{equation}}
\newcommand{\eeq}{\end{equation}}
\newcommand{\ba}{\begin{align}}
\newcommand{\ea}{\end{align}}

\renewcommand{\Re}{\mathop{\rm Re\,}}
\renewcommand{\Im}{\mathop{\rm Im\,}}

\DeclareMathOperator{\tr}{tr}

\DeclareMathOperator{\spec}{spec}

\DeclareMathOperator{\dom}{dom}

\DeclareMathOperator*{\nlim}{{\mathit n}-lim}
\DeclareMathOperator*{\slim}{{\mathit s}-lim}

\allowdisplaybreaks
\numberwithin{equation}{section}
\newtheorem{theorem}{Theorem}[section]
\newtheorem{lemma}[theorem]{Lemma}
\newtheorem{corollary}[theorem]{Corollary}
\newtheorem{hypothesis}[theorem]{Hypothesis}
\theoremstyle{definition}
\newtheorem{definition}[theorem]{Definition}
\newtheorem{example}[theorem]{Example}

\theoremstyle{remark}
\newtheorem{remark}[theorem]{Remark}

\begin{document}
\title[Monotonicity and Concavity]{Monotonicity and Concavity 
Properties of The Spectral Shift Function}
%\dedicatory{}
% Information for  author
\author[]{Fritz Gesztesy,
Konstantin A.~Makarov,  and  Alexander~K.~Motovilov  }
% Information for  author
%\author{Fritz Gesztesy}
\address{Department of Mathematics,
University of
Missouri, Columbia, MO
65211, USA}
\email{fritz@math.missouri.edu\newline
\indent{\it URL:}
http://www.math.missouri.edu/people/fgesztesy.html}
% Information for  author
%\author{Konstantin A.~Makarov}
\address{Department of Mathematics, University of
Missouri, Columbia, MO
65211, USA}
\email{makarov@azure.math.missouri.edu}
% Information for author
%\author{Alexander K.~Motovilov}
\address{
Physikalisches Institut, Universit\"at Bonn, D-53115~Bonn, Germany
 }
\email{motovilov@physik.uni-bonn.de}
\address{On leave of absence from the
 Laboratory of Theoretical Physics, JINR, 141980 Dubna, Russia
}
\email{motovilv@thsun1.jinr.ru}
%%%%%%%%%%%%%%%%%%%%%%%%%%%%%%%%%%%%%%%%%%%%%%%%%%%%%%%%%%%%%%%
\dedicatory{Dedicated with great pleasure to Sergio
Albeverio on the occasion of his 60th birthday}
%%%%%%%%%%%%%%%%%%%%%%%%%%%%%%%%%%%%%%%%%%%%%%%%%%%%%%%%%%%%%%%
%%%%%%%%%%%%%%%%
\date{September, 1999}
\subjclass{Primary 47B44, 47A10; Secondary 47A20, 47A40}

%%%%%%%%%%%%%%%%%%%%%%%%%%%%%%%%%%%%%%%%%%%%%%%%%%%%%%%%
\begin{abstract}
Let $H_0$ and $V(s)$ be self-adjoint, $V,V'$ continuously
differentiable in trace norm with $V''(s)\geq 0$ for
$s\in (s_1,s_2)$, and denote by
$\{E_{H(s)}(\lambda)\}_{\lambda\in\bbR}$ the family of spectral
projections of $H(s)=H_0+V(s)$. Then we prove for given
$\mu\in\bbR$, that $s\longmapsto
\tr\big (V'(s)E_{H(s)}((-\infty, \mu))\big ) $ is a
nonincreasing function with respect to $s$, extending a
result of Birman and Solomyak. Moreover, denoting by
$\zeta (\mu,s)=\int_{-\infty}^\mu d\lambda \,
\xi(\lambda,H_0,H(s))$ the integrated spectral shift
function for the pair $(H_0,H(s))$, we prove concavity of
$\zeta (\mu,s)$ with respect to $s$, extending previous 
results by
Geisler, Kostrykin, and Schrader. Our proofs employ 
operator-valued Herglotz functions and establish the latter 
as an effective tool in this context. 
\end{abstract}
%%%%%%%%%%%%%%%%%%%%%%%%%%%%%%%%%%%%%%%%%%%%%%%%%%%%%

\maketitle

%%%%%%%%%%%%%%%%%%%%%%%%%%%%%%%%%%%%%%%%%%%%%%%%%%%%%
\section{Introduction and principal results} \label{s1}
%%%%%%%%%%%%%%%%%%%%%%%%%%%%%%%%%%%%%%%%%%%%%%%%%%%%%%%

In the  following $\calH$ denotes a complex separable Hilbert
space with scalar product $(\,\cdot, \, \cdot)_{\calH}$ (linear
in the second factor) and norm $\|\cdot\|_\calH$, $\calB(\calH)$
represents the Banach space of bounded linear operators 
defined on
$\calH$,
$\calB_p(\calH), \,\, p\ge 1$ the standard Schatten-von Neumann
ideals of $\calB(\calH)$ (cf., e.g., \cite{GK69}, \cite{Si79})
 and $\bbC_+ $ (resp., $\bbC_-$) the open complex upper 
(resp., lower) 
half-plane. Moreover, real and imaginary parts of a bounded 
operator $T\in
\calB(\calH)$ are defined as usual by
$\Re (T)=(T+T^*)/2$, $\Im (T)=(T-T^*)/(2i)$.

The spectral shift function $\xi(\lambda, H_0,H)$ associated
with a pair of self-adjoint operators $(H_0,H)$, $H=H_0+V$,
$\dom(H_0)=\dom(H)$, where
\begin{equation}
\lb{nucl}
 V=V^*\in \calB_1(\calH),
\end{equation}
is one of the fundamental spectral characteristics in the
perturbation theory of self-adjoint operators.  It is 
well-known
(see \cite{Kr53}, \cite{Kr62}, \cite{Kr83}, \cite{Kr89},
\cite{Li52}) that
for a wide  function class ${\mathfrak K}(H_0,H)$, the 
Lifshits-Krein trace
formula holds, that is,
\begin{equation}
\label{5.1}
\tr(\varphi(H)-\varphi(H_0) )=\int_{\bbR}
d\lambda\, \varphi'(\lambda)\,\xi(\lambda, H_0,H)\,, 
\quad \varphi\in
{\mathfrak K}(H_0,H).
\end{equation}
In the case of trace class perturbations \eqref{nucl},
the spectral shift function is integrable, that is,
\begin{equation}\lb{L1}
\xi(\cdot, H_0, H)\in L^1(\bbR),
\end{equation}
and  the following relations hold
\begin{eqnarray}
\label{ksi1}
\|\xi(\cdot, H_0, H)\|_{L^1(\bbR)} 
&\le& \|V\|_{\calB_1(\calH)},\\
\label{trV}
\int_\bbR d\lambda\,\,\xi(\lambda, H_0, H)&=&\tr(V)\,.
\end{eqnarray}

The precise characterization of the class
\begin{equation}
{\mathfrak K}=\bigcap_{H_0, H}
{\mathfrak K}(H_0,H)  \lb{1.6}
\end{equation}
of all those $\varphi$ for which \eqref{5.1} holds for any pair
of self-adjoint operators $H_0$ and $H=H_0+V$ with a trace class
difference \eqref{nucl}, is still unknown. In particular, 
there are
functions $\varphi\in C^{1}_0(\bbR)$ for which \eqref{5.1} fails
(cf.~\cite{Fa80}, \cite{Pe94}). Necessary
conditions very close to sufficient ones for $\varphi$ belonging 
to the
class ${\mathfrak K}$ have been found by Peller \cite{Pe85},
\cite{Pe90}. Here we only note that \eqref{L1} and
$(\varphi(H)-\varphi(H_0))\in\calB_1(\calH)$ hold, and
\eqref{5.1}  is valid, if $\varphi'$  is the Fourier transform 
of a finite Borel measure,
\begin{equation}
\label{5.2}
\varphi'(\lambda)=\int_\bbR d \nu(t)\,e^{-it\lambda}, \quad
\varphi\in C^1(\bbR), \,\,\,\int_\bbR d |\nu(t)|<\infty.
\end{equation}
We denote the function class \eqref{5.2} by $\calW_1(\bbR)$.

Different representations for the spectral shift function  and 
their interrelationships  can be found in \cite{BP98}, for 
further information we refer to \cite[Ch.~19]{BW83}, 
\cite{BY93}, \cite{BY93a}, \cite{GMN99}, \cite[Ch.~8]{Ya92} and 
the references therein).

In the present short note we will focus on two particular 
results:  one, a monoto\-ni\-ci\-ty result obtained by Birman 
and Solomyak \cite{BS75}, the other, a concavity result obtained 
by Geisler, Kostrykin, and Schrader \cite{GKS95}, \cite{Ko99}.  
We also present some extensions and new proofs that we 
hope might give additional insights into the subject.

We start by recalling pertinent results discovered by Birman and 
Solomyak \cite{BS75} in connection with the spectral averaging 
formula (providing a representation for the spectral shift 
function via an integral over the coupling constant) and a 
monotonicity result of a certain trace with respect to the 
coupling constant parameter.

%%%%%%%%%%%%%%%%%%%%%%%%%%%%%%%%%%%%%%%%%%%%%%%%%%%%%%%%%%%%%%%
\begin{theorem}[\cite{BS75}]\label{t1.1}
Let $H_0$ and $V$ be self-adjoint in $\calH$,
$V\in\calB_1(\calH)$, and define
\begin{equation}
H_s=H_0+sV, \quad \dom(H_s)=\dom(H_0), \,\, s\in\bbR, \lb{1.8}
\end{equation}
with $\{E_{H_s}(\lambda)\}_{\lambda\in\bbR}$ the family of
orthogonal spectral projections of $H_s$. Moreover, denote by
$\xi(\cdot, H_0, H_1)$ the spectral shift function for the pair
$(H_0,H_1)$. Then for any Borel set $\Delta \subset \bbR$,
\begin{equation} \label{1.1}
\int_{\Delta} d\lambda\, \xi (\lambda , H_0, H_1)=
\int_0^1 ds\,\tr(VE_{H_s}(\Delta)).
\end{equation}
\end{theorem}
%%%%%%%%%%%%%%%%%%%%%%%%%%%%%%%%%%%%%%%%%%%%%%%%%%%%%%%%%%%%%%%

In the same paper \cite{BS75}, Birman and Solomyak proved another
remarkable statement concerning the monotonicity of the 
integrand in the
right-hand side  of (\ref{1.1}) with respect to $s$ for 
semi-infinite
intervals $\Delta =(-\infty,\lambda)$, $\lambda\in\bbR$.

%%%%%%%%%%%%%%%%%%%%%%%%%%%%%%%%%%%%%%%%%%%%%%%%%%%%%%%%%%%%%%%
\begin{theorem}[\cite{BS75}]\label{t1.2}
Assume the hypotheses in Theorem~\ref{t1.1}. Given $\mu\in 
\bbR$, the
function
\begin{equation}\label{1.2}
  s\longmapsto \tr\big (VE_{H_s}((-\infty,\mu))\big ), 
\quad s\in\bbR,
\end{equation}
is a nonincreasing function with respect to $s\in\bbR$.
\end{theorem}
%%%%%%%%%%%%%%%%%%%%%%%%%%%%%%%%%%%%%%%%%%%%%%%%%%%%%%%%%%%%%%%

The spectral averaging formula \eqref{1.1} combined with the
monotonicity result \eqref{1.2} is a convenient tool for 
producing
estimates for the spectral shift function \cite{BS75}. For 
instance,
\begin{equation}
\label{1.3}
\tr \big (VE_{H_1}((-\infty, \mu)) \big )\le
\int_{-\infty}^\mu d\lambda\, \xi (\lambda , H_0, H_1)
\le
\tr \big (VE_{H_0}((-\infty, \mu)) \big ).
\end{equation}
In particular, passing to the limit $\mu\to \infty$ in
\eqref{1.3} one obtains \eqref{trV} (see \cite{BS75}
for more details).

Another application of the pair of results
\eqref{1.1} and \eqref{1.2} leads to the proof of concavity
properties of the integrated spectral shift function with
respect to the coupling constant, originally discovered in the
case of Schr\"odinger operators by Geisler, Kostrykin, and
Schrader \cite{GKS95} and extended by Kostrykin \cite{Ko99}, 
\cite{Ko99a} to the general case presented next. 

%%%%%%%%%%%%%%%%%%%%%%%%%%%%%%%%%%%%%%%%%%%%%%%%%%%%%%%%%%%%%%
\begin{theorem}[\cite{Ko99}, \cite{Ko99a}] \label{t1.3}
Let $\xi(\cdot, H_0, H_s)$ be the spectral shift function
in Theorem~\ref{t1.1}. Given $\mu\in \bbR$, the integrated
spectral shift function
\begin{equation}
\zeta_s(\mu)=
\int_{-\infty}^\mu d\lambda\, \xi (\lambda , H_0, H_s),
\quad s\in\bbR
\end{equation}
is a concave function with respect to the coupling constant 
$s\in\bbR$.
More precisely, for any $s, t \in \bbR$, and for all 
$\alpha\in[0,1]$,
the following inequality
\begin{equation}
\zeta_{\alpha s+ (1-\alpha)t}(\mu)\ge \alpha
\,\zeta_s(\mu) +(1-\alpha) \,\zeta_t(\mu)
\end{equation}
holds. Moreover, $\zeta_s(\mu)$ is subadditive with respect to
 $s\in (0,\infty)$ in the sense that for any $s, t\ge 0$,
\begin{equation}
\zeta_{s+t}(\mu)\le\zeta_s(\mu)+\zeta_t(\mu).
\end{equation}
\end{theorem}
%%%%%%%%%%%%%%%%%%%%%%%%%%%%%%%%%%%%%%%%%%%%%%%%%%%%%%%%%%%%%%%%

While Theorem~\ref{t1.3} focuses on a linear coupling constant 
dependence in $H_s=H_0+sV$, Kostrykin~\cite{Ko99} also 
discusses the case of a nonlinear dependence on $s$ for 
operators of the form $H(s)=H_0+  V(s)$:

%%%%%%%%%%%%%%%%%%%%%%%%%%%%%%%%%%%%%%%%%%%%%%%%%%%%%%%%%%%%%%%%
\begin{theorem}[\cite{Ko99}, \cite{Ko99a}] \label{t1.3a}
Suppose $f:\bbR\to\bbR$ is an nonincreasing function of bounded 
variation and $\{V(s)\}_{s\in\bbR}\in\calB_1(\calH)$ is operator 
concave (i.e., $V(\alpha s+(1-\alpha)t)\geq \alpha V(s)
+(1-\alpha)V(t)$ for all $\alpha\in [0,1]$, $s,t\in\bbR$).  Then 
\begin{equation}
s \longmapsto g(V(s))=\int_\bbR d\lambda\,f(\lambda)
\xi(\lambda,H_0,H_0+V(s)), 
\end{equation}
is concave in $s\in\bbR$. More precisely, for
all $0\le\alpha\le1$ and all $s,t\in\bbR$, the following 
inequality
\begin{equation} \label{Concavity}
g(V(\alpha s+(1-\alpha)t))\geq \alpha g(V(s))+(1-\alpha)g(V(t))
\end{equation}
holds.
\end{theorem}
%%%%%%%%%%%%%%%%%%%%%%%%%%%%%%%%%%%%%%%%%%%%%%%%%%%%%%%%%%%%%%

Actually, Kostrykin considered the general case of relative 
trace class 
perturbations in \cite{Ko99} but we omit further details in 
this note.

%%%%%%%%%%%%%%%%%%%%%%%%%%%%%%%%%%%%%%%%%%%%%%%%%%%%%%%%%%%%%%%
\begin{remark} \lb{r1.4}
The results of  Theorems~\ref{t1.1} and \ref{t1.2} in \cite{BS75}
have
been obtained using the approach of Stieltjes' double operator
integrals \cite{BS66}--\cite{BS73}.  Birman and Solomyak treated
the case $V(s)=sV$, $V\in\calB_1(\calH)$,  that is, they 
discussed
the case of a
linear dependence of the perturbation $V(s)$ with respect to the
coupling constant parameter $s$.  The general case of a
nonlinear dependence $V(s)$ of $s$, assuming
$V'(s)\ge0$, in the context of the spectral averaging result
\eqref{1.1} has recently been treated by Simon \cite{Si98}.
In Theorem~\ref{t1.6} below we cite the most recent result of
this type obtained in \cite{GMN99}.
\end{remark}
%%%%%%%%%%%%%%%%%%%%%%%%%%%%%%%%%%%%%%%%%%%%%%%%%%%%%%%%%%%%%%

It is convenient to introduce the following hypothesis.

%%%%%%%%%%%%%%%%%%%%%%%%%%%%%%%%%%%%%%%%%%%%%%%%%%%%%%%%%%%%%
\begin{hypothesis}\label{h1.5}
Let $H_0$ be a self-adjoint operator in $\calH$ with domain
$\dom(H_0)$, and assume $\{V(s)\}_{s\in \Omega}\subset
\calB_1(\calH)$ to be a family of self-adjoint trace class
operators in $\calH$, where $\Omega \subseteq \bbR$ denotes
an open interval with $0\in\Omega$. Moreover, suppose that 
$V(s)$ is
continuously differentiable in
$\calB_1(\calH)$-norm with respect to $s\in\Omega$.
For  convenience (and without loss of generality) we may assume
that $V(0)=0$ in the following.
\end{hypothesis}
%%%%%%%%%%%%%%%%%%%%%%%%%%%%%%%%%%%%%%%%%%%%%%%%%%%%%%%%%%%%

In the rest of the paper we will frequently use the notation 
$(s_1,s_2)\subset\subset \Omega$ to denote an open interval 
that is strictly contained in the interval $\Omega=(a,b)$ 
(i.e., $a<s_1<s_2<b$).

%%%%%%%%%%%%%%%%%%%%%%%%%%%%%%%%%%%%%%%%%%%%%%%%%%%%%
\begin{theorem} [\cite{GMN99}] \label{t1.6}
Assume Hypothesis~\ref{h1.5} and $0\in (s_1,s_2)\subset
\subset \Omega$. Let
\begin{equation}
H(s)=H_0+V(s), \quad \dom(H(s))=\dom(H_0), \quad s\in (s_1,s_2),
\end{equation}
with $\{E_{H(s)}(\lambda)\}_{\lambda\in\bbR}$ the family of
orthogonal spectral projections of $H(s)$ and denote by
$\xi(\cdot,H_0,H(s))$ the spectral shift function for the pair
$(H_0,H(s))$.  Then for any Borel set $\Delta\subset\bbR$  the
following spectral averaging formula holds
\begin{equation}
\int_{\Delta} d\lambda\, \big (\xi (\lambda , H_0, H(s_2))-
\xi (\lambda , H_0, H(s_1))\big )=
\int_{s_1}^{s_2} ds\,\tr(V'(s)E_{H(s)}(\Delta)).
\end{equation}
\end{theorem}
%%%%%%%%%%%%%%%%%%%%%%%%%%%%%%%%%%%%%%%%%%%%%%%%%%%%%%%%%

The principal new result of the present note is an
extension of the monotonicity result, Theorem~\ref{t1.2},
to the case of a nonlinear dependence of $V(s)$ on
$s$. In particular, we provide a new strategy of proof for 
such results, which appears to be interesting in itself.

%%%%%%%%%%%%%%%%%%%%%%%%%%%%%%%%%%%%%%%%%%%%%%%%%%%%%%%%%
\begin{theorem} \label{t1.7}
Assume Hypothesis~\ref{h1.5} and $0\in (s_1,s_2)\subset
\subset \Omega$.
Suppose in addition, that the derivative  $V'(s)$ is
continuously differentiable in
$\calB_1(\calH)$-norm with respect to $s\in(s_1,s_2)$
and that $V(s)$ is concave in the sense that
\begin{equation}
     0\ge V''(s) \in \calB_1(\calH), \quad s\in (s_1, s_2).
\end{equation}
Then, given $\mu\in \bbR$, the function
\begin{equation}
s\longmapsto \tr\big (V'(s)E_{H(s)}((-\infty, \mu))\big ),
\quad s\in (s_1, s_2),
\end{equation}
is a nonincreasing function with respect to $s\in (s_1, s_2)$.
\end{theorem}
%%%%%%%%%%%%%%%%%%%%%%%%%%%%%%%%%%%%%%%%%%%%%%%%%%%%%%%%%

Combining Theorems~\ref{t1.6} and~\ref{t1.7} one obtains the 
following result.

%%%%%%%%%%%%%%%%%%%%%%%%%%%%%%%%%%%%%%%%%%%%%%%%%%%%%%%%%%%%%
\begin{corollary} \label{c1.8}
Suppose the hypotheses of Theorem~\ref{t1.7}.
Then, for given $\mu\in \bbR$, the integrated  spectral shift
function
\begin{equation}
\label{ssf}
\zeta(\mu,s)=
\int_{-\infty}^\mu d\lambda\, \xi(\lambda, H_0, H(s))
\end{equation}
is  concave in $s\in (s_1,s_2)$. More precisely, for
all $0\le\alpha\le1$ and all
$s,t\in (s_1,s_2)$, the following inequality
\begin{equation}
\label{Concav}
\zeta(\mu,\alpha s+ (1-\alpha)t )\ge \alpha
\,\zeta(\mu,s) +(1-\alpha) \,\zeta(\mu,t)
\end{equation}
holds. Moreover, $\zeta(\mu,s)$ is subadditive with respect
to $s\in [0,s_2)$ in the sense that for any $s, t\ge 0$, 
$s+t\in [0,s_2)$,
\begin{equation}
\label{3.16a}
\zeta(\mu, s+t)\le\zeta(\mu,s)+
\zeta(\mu, t).
\end{equation}
\end{corollary}
%%%%%%%%%%%%%%%%%%%%%%%%%%%%%%%%%%%%%%%%%%%%%%%%%%%%%%%%%%%

We emphasize that Corollary~\ref{c1.8} is a special case 
of Kostrykin's Theorem~\ref{t1.3a}.

\medskip

As explained in \cite{Bi94}, \cite{BS75}, and \cite{BY93},
the original proofs of Theorems~\ref{t1.1} and~\ref{t1.2}
in \cite{BS75} were motivated by a real analysis approach
to the spectral shift function in contrast to M.~Krein's
complex analytic treatment. In this note we return to complex
analytic proofs in the  spirit of M.~Krein and provide a
proof of the monotonicity result, Theorem~\ref{t1.7}, based
on operator-valued Herglotz function techniques. For
various recent applications of  this formalism we refer
to \cite{GKMT98}, \cite{GM99}, \cite{GM99a}, \cite{GMN99},
and \cite{GT97}.

%%%%%%%%%%%%%%%%%%%%%%%%%%%%%%%%%%%%%%%%%%%%%%%%%%%%%%%%%
%%%%%%%%%%%%%%%%%%%%%%%%%%%%%%%%%%%%%%%%%%%%%%%%%%%%%%%%%
\section{A property of operator-valued Herglotz functions}
%%%%%%%%%%%%%%%%%%%%%%%%%%%%%%%%%%%%%%%%%%%%%%%%%%%%%%%%%
%%%%%%%%%%%%%%%%%%%%%%%%%%%%%%%%%%%%%%%%%%%%%%%%%%%%%%%%%

We recall that $f:\bbC_+\to\bbC$ is called a {\it Herglotz 
function} if it is analytic and $f(\bbC_+)\subseteq\bbC_+$. 
In this case we extend $f$ to $\bbC_-$ in the usual manner, 
that is, by $f(\overline z)=\overline{f(z)}$, $z\in\bbC_+$.

The principal purpose of this section is to obtain some 
generalizations of the following elementary result.

%%%%%%%%%%%%%%%%%%%%%%%%%%%%%%%%%%%%%%%%%%%%%%%%%%%%%%%%%
\begin{lemma} \lb{l2.1}
Let $P$ and $Q$ be two rational Herglotz functions vanishing at 
infinity  and let $\Gamma$ be a closed clockwise oriented Jordan 
contour encircling some of the poles of $P$ and $Q$ starting 
from the left (and without any poles of $P$ and $Q$ on 
$\Gamma$). Then
\begin{equation}
\frac{1}{2\pi i} \oint_{\Gamma} dz \,
P(z)Q(z)\ge 0.
\end{equation}
\end{lemma}
%%%%%%%%%%%%%%%%%%%%%%%%%%%%%%%%%%%%%%%%%%%%%%%%%%%%%%%%%
\begin{proof}
By the hypotheses on $P$ and $Q$ we may write
\begin{align}
P(z)&=\sum_{j\in J_1}A_j(p_j-z)^{-1}, \quad A_j\geq 0,
\,\, p_j\in\bbR, \, j\in J_1, \\
Q(z)&=\sum_{\ell\in J_2}B_\ell (q_\ell -z)^{-1},
\quad B_\ell\geq 0, \,\, q_\ell\in\bbR, \, \ell\in J_2, \\
\end{align}
with $J_1,J_2$ finite index sets. Next one decomposes $P$
and $Q$ with respect to their poles located in the interior
and exterior of the bounded domain encircled by $\Gamma$,
\begin{align}
P(z)&=P_{\text{int}}(z)+P_{\text{ext}}(z) \no \\
&=\sum_{j\in J_{1,\text{int}}}A_j(p_j-z)^{-1}
+\sum_{j\in J_{1,\text{ext}}}A_j(p_j-z)^{-1}, \\
Q(z)&=Q_{\text{int}}(z)+Q_{\text{ext}}(z) \no \\
&=\sum_{\ell\in J_{2,\text{int}}}B_\ell (q_\ell -z)^{-1}
+\sum_{\ell\in J_{2,\text{ext}}}B_\ell (q_\ell -z)^{-1}, \\
&\hspace*{2.7cm} J_k=J_{k,\text{int}}\cup J_{k,\text{ext}},
\,\, k=1,2.
\end{align}
Then straightforward residue computations yield
\begin{align}
&\frac{1}{2\pi i} \oint_{\Gamma} dz \,P(z)Q(z) \no \\
&=\frac{1}{2\pi i} \oint_{\Gamma} dz \,
P_{\text{int}}(z)Q_{\text{int}}(z)
+\frac{1}{2\pi i} \oint_{\Gamma} dz \,
P_{\text{int}}(z)Q_{\text{ext}}(z) \no \\
&\quad +\frac{1}{2\pi i} \oint_{\Gamma} dz \,
P_{\text{ext}}(z)Q_{\text{int}}(z)
+\frac{1}{2\pi i} \oint_{\Gamma} dz \,
P_{\text{ext}}(z)Q_{\text{ext}}(z) \lb{2.8} \\
&=\frac{1}{2\pi i} \oint_{\Gamma} dz \,
\sum_{j\in J_{1,\text{int}}}
\sum_{\ell\in J_{2,\text{int}}}
A_jB_\ell (p_j-z)^{-1}(q_\ell-z)^{-1} \no \\
&\quad +\sum_{j\in J_{1,\text{int}}} A_jQ_{\text{ext}}(p_j)
+\sum_{\ell\in J_{2,\text{int}}} B_\ell P_{\text{ext}}
(q_\ell) \no \\
&=\frac{1}{2\pi i} \oint_{\Gamma} dz \,
\sideset{}{'}\sum_{j\in J_{1,\text{int}}, 
\ell\in J_{2,\text{int}}}
A_jB_\ell (q_\ell -p_j)^{-1}\big((q_\ell-z)^{-1}
-(p_j-z)^{-1}\big) \no \\
&\quad +\sum_{j\in J_{1,\text{int}}} A_jQ_{\text{ext}}(p_j)
+\sum_{\ell\in J_{2,\text{int}}} B_\ell P_{\text{ext}}
(q_\ell) \lb{2.9} \\
&=\sum_{j\in J_{1,\text{int}}} A_jQ_{\text{ext}}(p_j)
+\sum_{\ell\in J_{2,\text{int}}} B_\ell P_{\text{ext}}
(q_\ell) \no \\
&=\sum_{j\in J_{1,\text{int}}}\sum_{\ell\in J_{2,\text{ext}}}
A_jB_\ell(q_\ell-p_j)^{-1}
+\sum_{j\in J_{1,\text{ext}}}\sum_{\ell\in J_{2,\text{int}}}
A_jB_\ell (p_j-q_\ell)^{-1}\geq 0.
\end{align}
Here the last integral in \eqref{2.8} vanishes since the
integrand is analytic inside $\Gamma$ and the first integral
in \eqref{2.9} vanishes by symmetry. Moreover, we used the
symbol $\sideset{}{'}\sum$ to indicate summation only over
those $j$ and $\ell$ with $p_j\neq q_\ell$, since only
first-order poles contribute in this calculation.
\end{proof}
%%%%%%%%%%%%%%%%%%%%%%%%%%%%%%%%%%%%%%%%%%%%%%%%%%%%%%%%%%%%%%

Next we turn to operator-valued extensions of the concept
of Herglotz functions.

%%%%%%%%%%%%%%%%%%%%%%%%%%%%%%%%%%%%%%%%%%%%%%%%%%%%%%%%%%%%%%
\begin{definition} \label{d2.2}
$M:\bbC_+\to \calB(\calH)$ is called an {\it
operator-valued Herglotz function} if $M$ is analytic on
$\bbC_+$ and $\Im (M(z))\ge 0$ for all $z\in \bbC_+$.
\end{definition}
%%%%%%%%%%%%%%%%%%%%%%%%%%%%%%%%%%%%%%%%%%%%%%%%%%%%%%%%%

Any operator-valued Herglotz function admits a canonical
representation, which can be considered a generalization of
the dilation theory of maximal dissipative operators to the case
of operator-valued Herglotz functions. In the following, 
however, we will
focus on Herglotz functions of the resolvent-type
\begin{equation}
M(z)=K^*(L-z)^{-1}K,\quad z\in \bbC_+,
\end{equation}
where $K$ is a bounded operator between the Hilbert spaces 
$\calK$
and $\calH$, $K\in\calB(\calH,\calK)$, $\calK\supseteq\calH$,
while  $L$ is assumed to be a
self-adjoint operator in $\calK$, bounded from below, and
with a gap in its spectrum as described in Theorem~\ref{t2.3} 
below.

%%%%%%%%%%%%%%%%%%%%%%%%%%%%%%%%%%%%%%%%%%%%%%%%%%%%%%%%%%%%%%%
\begin{theorem} \label{t2.3}
Let $M_j:\bbC_+\to \calB_2(\calH)$, $j=1,2$, be  operator-valued
Herglotz functions (taking values in the space of the
Hilbert-Schmidt operators) admitting the representations
\begin{equation}
\label{0.7}
M_j(z)=K_j(L_j-z)^{-1}K_j^*,
\end{equation}
where $\calK_j$ and $\calH$, $\calK_j\supseteq\calH$, are
Hilbert spaces, $L_j$ are self-adjoint operators in $\calK_j$
bounded from below, and $K_j\in \calB_4(\calK_j,\calH)$,
$j=1,2$. Suppose $\calD$ is a domain in the complex
plane and  $ (a,b)$ an
open interval such that
\begin{equation}
\label{0.5}
     a<\min_{j=1,2} \inf\bigl(\spec(L_j)\bigr)
\end{equation}
and
\begin{equation}
\label{0.6}
\bigl\{\spec(L_1)\cup\spec(L_2)\bigr\}\cap\calD\subset(a,b)
\subset\calD.
\end{equation}
In addition, assume that $\Gamma$ is a closed oriented
Jordan contour in
$\calD$ encircling the interval $[a,b]$ in the clockwise 
direction, and
$\varphi$   an analytic function on $\calD$, nonnegative and
nonincreasing on $(a,b)$, that is,
\begin{equation}
\label{0.8}
\varphi\bigl|_{(a,b)}\bigr.\ge 0,
\end{equation}
and
\begin{equation}
\varphi'\bigl|_{(a,b)}\bigr.\le 0.
\end{equation}

Then
\begin{equation}
\label{0.10}
\frac{1}{2\pi i} \oint_{\Gamma} dz \,\varphi
(z)\,\tr\big (M_1(z)M_2(z) \big )\ge 0.
\end{equation}
\end{theorem}
%%%%%%%%%%%%%%%%%%%%%%%%%%%%%%%%%%%%%%%%%%%%%%%%%%%%%%%%%
\begin{proof}
Given $n\in \bbN$, introduce the partition $a=t_0<t_1<t_2 <
\ldots< t_{n-1}<t_n=b$ of the closed interval $[a,b]$,
\begin{equation}
t_k=a+k\frac{b-a}{n}, \quad k=1, 2, \,\ldots\, , n,
\end{equation}
and denote by $\chi^{(n)}(\lambda)$ the piecewise continuous
function
\begin{equation}
\chi^{(n)}(\lambda)=\sum_{k=1}^n \,\,t_k\,\,
\chi_{[t_{k-1},t_k)}(\lambda), \quad \lambda\in \bbR,
\end{equation}
where  $\chi_\Delta(\cdot)$ is the characteristic function
of the set $\Delta\subset\bbR$.

In the Hilbert space $\calK_j$ introduce the (possibly
unbounded) operators
\begin{align}
&L_j^{(n)}=\int_{(a,b)}\chi^{(n)}(\lambda)\,\,dE_{L_j}(\lambda)
 +\int_{(b,\infty)}\lambda\,\,dE_{L_j}(\lambda),
\label{0.13} \\
&\dom(L_j^{(n)})=\dom(L_j), \quad  n\in \bbN, \,\, j=1,2.
\no
\end{align}
We note, that by definition \eqref{0.13}, the spectrum of $L_j$,
$j=1,2$, in the interval $(a,b)$ consists of finitely many
eigenvalues (possibly of infinite multiplicity).

By \eqref{0.6} one infers that the sequence of operators
$\{L_j^{(n)}\}_{n=1}^\infty $ converges in norm resolvent
sense to  $L_j $ in the Hilbert space $\calK_j$, $j=1,2$. This
convergence, in turn, combined with the hypothesis $K_j\in
\calB_4(\calK_j,\calH)$, $j=1,2$,  implies the convergence
\begin{equation}
\lim_{n\to \infty} \frac{1}{2\pi i} \oint_{\Gamma} dz
\,\varphi(z)\,\tr\big (M_1^{(n)}(z)M_2^{(n)}(z) \big )
=\frac{1}{2\pi i} \oint_{\Gamma} dz \,\varphi(z)
\,\tr\big(M_1(z)M_2(z)\big), \label{0.14}
\end{equation}
where in obvious notation
\begin{equation}
M_j^{(n)}(z)=K_j (L_j^{(n)}-z)^{-1}K_j^*, \quad j=1,2.
\end{equation}
Thus, in order to prove \eqref{0.10} it suffices to check that
every term on the left-hand side of \eqref{0.14} is nonnegative.

In the following it is useful to decompose the Herglotz functions
$M_j^{(n)}(z)$ in the form
\begin{equation}
\label{0.16}
   M_j^{(n)}(z)=N_j^{(n)}(z)+\widetilde N_j^{(n)}(z), \quad j=1,2,
\end{equation}
where
\begin{align}
 N_j^{(n)}(z)&= K_j E_{L_j}\,\bigl((a,b)\bigr)\,
(L_j^{(n)}-z)^{-1}E_{L_j}\bigl((a,b)\bigr)\,K_j^*, \\
\widetilde{N}_j^{(n)}(z)&=K_j\,E_{L_j}\bigl([b, \infty)\bigr)\,
(L_j^{(n)}-z)^{-1}
E_{L_j}\bigl([b,\infty)\bigr)\,K_j^*, \quad j=1,2 \label{0.18}
\end{align}
are Herglotz functions associated with the spectral
subspaces $E_{L_j}\bigl((a,b)\bigr)\,\calK_j$ and
$E_{L_j}\bigl([b,\infty)\bigr)\,\calK_j$ of $L_j$, $j=1,2$.

According to the decomposition \eqref{0.13}, the Herglotz
functions $N_j^{(n)}(z)$ are rational operator-valued
functions of the form
\begin{align}
\label{0.19}
N_j^{(n)}(z)=
\sum_{k=1}^n\,\,\frac {Q_j^{(n), k}}{\,\,t_k-z\,\,}, 
\quad j=1,2,
\end{align}
where
\begin{equation}
\label{0.20}
    Q_j^{(n), k}=K_j E_{L_j}\bigl([t_{k-1},t_k)\bigr)\,
    K_j^*\ge 0, \quad
    k=1,\, \ldots\, , n, \,\, j=1,2.
\end{equation}
Given $n\in\bbN$, one decomposes the integrals on the
left-hand side of
\eqref{0.14} as a sum of four terms
\begin{eqnarray}
\frac{1}{2\pi i} \oint_{\Gamma} dz \,\varphi
(z)\,\tr\big(M_1^{(n)}(z)M_2^{(n)}(z)\big)
&=&\frac{1}{2\pi i} \oint_{\Gamma} dz \,\varphi(z)\,\tr
\big(N_1^{(n)}(z)N_2^{(n)}(z)\big) \no \\
&& +\frac{1}{2\pi i} \oint_{\Gamma} dz \,\varphi(z)\,\tr
\big(\widetilde N_1^{(n)}(z) \widetilde N_2^{(n)}(z)\big) \no \\
&& +\frac{1}{2\pi i} \oint_{\Gamma} dz \,\varphi(z)
\,\tr\big ( \widetilde N_1^{(n)}(z)N_2^{(n)}(z)\big ) \no \\
&& +\frac{1}{2\pi i} \oint_{\Gamma} dz \,\varphi(z)
\,\tr\big (N_1^{(n)}(z) \widetilde N_2^{(n)}(z)\big) \no \\
&=& J_1+J_2+J_3+J_4. \label{0.21}
\end{eqnarray}
According to \eqref{0.19}, the first integral $J_1$  in
\eqref{0.21}  can be represented as follows
\begin{align}
J_1&=\frac{1}{2\pi i}\sum_{k, m}^n \oint_{\Gamma} dz \,
\frac{\varphi (z)}{(t_k-z)(t_m-z)} \,
\tr(Q_1^{(n), k}Q_2^{(n),m})\no \\
&=\frac{1}{2\pi i}
\sum^n_{\mbox{\scriptsize$\begin{array}{c} 
k,m\\k\neq m\end{array}$}}
\oint_{\Gamma} dz \,
\frac{\varphi (z)}{(t_k-z)(t_m-z)} \,
\tr(Q_1^{(n), k}Q_2^{(n),m})  \no \\
&\quad +\frac{1}{2\pi i}\sum_{k=1}^n \oint_{\Gamma} dz \, 
\frac{\varphi
(z)}{(t_k-z)^2} \,\tr(Q_1^{(n), k}Q_2^{(n),k}). \label{0.22}
\end{align}
Applying the residue theorem  to the integrals on
the right-hand side of (\ref{0.22}) one infers
\begin{align}
J_1&=
\sum^n_{\mbox{\scriptsize$\begin{array}{c} 
k,m\\k\neq m\end{array}$}}
\bigg(-\frac{\varphi(t_k)-\varphi(t_m)}{t_k-t_m}\bigg)
\,\tr(Q_1^{(n), k}Q_2^{(n),m}) \no \\
&\quad -\sum^n_{k=1}
\,\varphi'(t_k)
\,\tr(Q_1^{(n), k}Q_2^{(n),k})\ge 0, \label{2.222}
\end{align}
since  $\varphi$ is a nonincreasing differentiable function 
and since
the inequalities
\begin{equation}
\tr(Q_1^{(n),k}Q_2^{(n),m})\ge0,\qquad k,m=1,\,\ldots\,,n
\end{equation}
hold (due to the fact that the operators $Q_j^{(n), k}$,
 $k=1,2,\ldots, n$, $j=1,2$, are nonnegative by \eqref{0.20}).

The  integral $J_2$ vanishes,
\begin{equation}
\label{2.28}
J_2=0,
\end{equation}
since
$\tr\big(\widetilde{N}_1^{(n)}(z)\widetilde N_2^{(n)}(z) \big )$
is holomorphic in $\calD$.

The remaining integrals $J_3$ and $J_4$ can also be evaluated by
the residue theorem and one obtains
\begin{align}
J_3&=\sum_{k=1}^n\,\varphi(t_k)\,\tr( \widetilde
N_1^{(n)}(t_k)Q_2^{(n),k}), \label{0.29} \\
J_4&=\sum_{k=1}^n\,\varphi(t_k)\,\tr(Q_1^{(n),k}\widetilde
N_2^{(n)}(t_k)). \label{0.30}
\end{align}
We recall that $\varphi(t_k)\ge 0$, $k=1,\ldots,n$, by \eqref{0.8}.
At the same time
\begin{equation}
\widetilde{N}_j^{(n)}(t_k)\ge 0, \quad k=1,2,\,
\ldots\,, n, \,\, n\in \bbN, \,\, j=1,2,
\end{equation}
by \eqref{0.13} and \eqref{0.18}, while
\begin{equation}
{Q}_j^{(n),k}\ge 0, \quad k=1,2,\,
\ldots\,, n, \,\, n\in \bbN, \,\, j=1,2,
\end{equation}
by \eqref{0.20}. Hence,
\begin{equation}
\tr(\widetilde{N}_1^{(n)}(t_k)Q_2^{(n),k})\ge0
\text{ and }
\tr(Q_1^{(n),k}\widetilde{N}_2^{(n)}(t_k))\ge0,
\quad k=1,\,\ldots\,,n,
\end{equation}
and thus, $J_3\geq0$ and $J_4\geq0$. Together with
\eqref{2.222} and \eqref{2.28} (combined with \eqref{0.14})
this proves \eqref{0.10}.
\end{proof}
%%%%%%%%%%%%%%%%%%%%%%%%%%%%%%%%%%%%%%%%%%%%%%%%%%%%%%%%%%%%%%%
\begin{remark}\label{r2.4}
(i) Conditions $K_j\in{\mathcal B}_4(\calK_j, \mathcal H)$,
$j=1,2$, in Theorem~\ref{t2.3} can be relaxed.  In fact,  it
suffices  to require that
\begin{equation}
M_1(z)M_2(z)\in \calB_1(\calH), \quad z\in \bbC_+.
\end{equation}

%%%%%%%%%%%%%%%%%%%%%%%%%%%%%%%%%%%%%%%%%%%%%%%%%%%%%%%%%%%%%%%
\smallskip

\noindent (ii)
Hypotheses \eqref{0.5} and \eqref{0.6} concerning the
semiboundedness of $L_j$, $j=1,2$, and the existence of a gap in
their spectra in Theorem~\ref{t2.3}, are also unnecessarily
stringent. In
fact, it  suffices to assume  that
the holomorphy domains of the Herglotz functions $M_j(z)$
$j=1,2$ given by \eqref{0.7} include the set
$(-\infty,\alpha)\cup(\beta,\gamma)$ for some $a<\alpha$ and
$b<\beta<\gamma$.

%%%%%%%%%%%%%%%%%%%%%%%%%%%%%%%%%%%%%%%%%%%%%%%%%%%%%%%%%%
\smallskip

\noindent (iii) If the operators $L_j$ are bounded and
\begin{equation}
\label{0.33}
\spec(L_j)\subset (a,b),  \quad j=1,2,
\end{equation}
then the terms \eqref{0.18} in the decomposition \eqref{0.16}
vanish and hence the integrals \eqref{0.29} and
\eqref{0.30} vanish too. This means that  under assumption
\eqref{0.33}, condition \eqref{0.8} is redundant.
\end{remark}
%%%%%%%%%%%%%%%%%%%%%%%%%%%%%%%%%%%%%%%%%%%%%%%%%%%%%%%%%%%%%%

%%%%%%%%%%%%%%%%%%%%%%%%%%%%%%%%%%%%%%%%%%%%%%%%%%%%%%%%%%%%%%%
\section{ Monotonicity,  Concavity, and Subadditivity}
%%%%%%%%%%%%%%%%%%%%%%%%%%%%%%%%%%%%%%%%%%%%%%%%%%%%%%%%%%%%%%%

Throughout this section we assume  Hypothesis~\ref{h1.5} and
 recall that
\begin{equation} \label{Hs}
H(s)=H_0+V(s), \quad \dom(H(s))=\dom(H_0), \,\,
s\in\Omega.
\end{equation}

First, we  treat the case of  bounded  $H_0$,
$H_0\in\calB(\calH)$, and
study differential and monotonicity properties of the function
\begin{equation}
    s\longmapsto\tr\bigl(V'(s)\,\varphi(H(s))\bigr),
    \quad  s\in  \Omega,
\end{equation}
where $\varphi$ is analytic on a domain $\calD$ containing
the spectra of the family \eqref{Hs}
\begin{equation}
  \bigcup_{s\in \Omega}\bigl(\spec(H(s))\bigr)\subset\calD\,.
\end{equation}

Next we introduce the following additional hypothesis, which is
motivated in part by Remark~\ref{r2.4}\,(iii).

%%%%%%%%%%%%%%%%%%%%%%%%%%%%%%%%%%%%%%%%%%%%%%%%%%%%%%%%%%%%%%%%%%
\begin{hypothesis}\label{h3.1}
Let $H$ be a bounded self-adjoint operator, $\calD$ a domain 
of the
complex plane, $(a,b)\subset\bbR$ an open interval, and $\Gamma$ a
closed clockwise oriented Jordan contour in $\calD$ encircling the
interval $[a,b]$ such that
\begin{equation}
\spec(H)\subset(a,b)\subset\calD.
\end{equation}
\end{hypothesis}
%%%%%%%%%%%%%%%%%%%%%%%%%%%%%%%%%%%%%%%%%%%%%%%%%%%%%%%%%%%%%%%

The following remark shows that  Hypothesis~\ref{h3.1} is 
stable under
 small (compact) perturbations of $H$.

%%%%%%%%%%%%%%%%%%%%%%%%%%%%%%%%%%%%%%%%%%%%%%%%%%%%%%%%%%%%%%%
\begin{remark}\label{r3.2}
Suppose that the collection $\{H(s_0),\calD,(a,b),\Gamma\}$ 
satisfies
Hypothesis~\ref{h3.1} for some
$s_0\in(s_1,s_2)$. By perturbation arguments, one infers the
existence of a neighborhood $S$  of $s_0$ such that
$\{H(s),\calD,(a,b),\Gamma\}$
 also satisfy Hypothesis~\ref{h3.1}
for $s\in S$.  Moreover, if $\varphi$ is  analytic on $\calD$,
then the (bounded) operators $\varphi(H(s))$, $s\in S$, are
well-defined by the Riesz integral
\begin{equation}
\label{2.18}
 \varphi(H(s))= \frac{1}{2\pi i} \oint_{\Gamma} dz \,\varphi(z)
 (H(s)-z)^{-1}, \quad s\in S.
\end{equation}
\end{remark}
%%%%%%%%%%%%%%%%%%%%%%%%%%%%%%%%%%%%%%%%%%%%%%%%%%%%%%%%%%%%%%%%

%%%%%%%%%%%%%%%%%%%%%%%%%%%%%%%%%%%%%%%%%%%%%%%%%%%%%%%%%%%%%%%%
\begin{lemma}\label{l3.3}
Assume Hypothesis~\ref{h1.5} and  let $S$ be the
neighborhood of $s_0$ in Remark~\ref{r3.2}. Suppose in
addition that
$V'(s)$ is continuously differentiable in
$\calB_1(\calH)$-norm for $s\in S$ and
that the collection $\{H(s_0),\calD,(a,b),\Gamma\}$
satisfies Hypothesis~\ref{h3.1} for some
$s_0\in(s_1,s_2)$.
If $\varphi$ is an analytic function on $\calD$, then
\begin{align}
\label{2.19}
    s\longmapsto \tr\bigl(V'(s)\,\varphi(H(s))\bigr), 
\quad s\in S,
\end{align}
is differentiable on $S$  and
\begin{eqnarray}
\nonumber
\frac{d}{ds}\tr\bigl(V'(s)\,\varphi(H(s))\bigr)&=&
 \tr\bigl(V''(s)\,\varphi(H(s))\bigr)\\
\label{2.20}
&&-\frac{1}{2\pi i}
\oint_{\Gamma}dz\,\varphi(z)\,\tr\bigl[V'(s)(H(s)-z)^{-1}\bigr]^2.
\end{eqnarray}
\end{lemma}
%%%%%%%%%%%%%%%%%%%%%%%%%%%%%%%%%%%%%%%%%%%%%%%%%%%%%%%%%%%%%%%%
\begin{proof}
By Remark~\ref{r3.2}, the operators  $\varphi(H(s))$, $s\in S$, 
are
well-defined (cf.~\eqref{2.18}). In particular,
\begin{align}
&\frac{\varphi(H(s))-\varphi(H(s_0))}{s-s_0} =
\frac{1}{2\pi i} \oint_{\Gamma} dz \,\varphi
(z)\frac{(H(s)-z)^{-1}-(H(s_0)-z)^{-1}}{s-s_0} \no \\
&=-\frac{1}{2\pi i} \oint_{\Gamma} dz \,\varphi
(z)(H(s)-z)^{-1}\frac{V(s)-V(s_0)}{s-s_0}(H(s_0)-z)^{-1}\,.
\label{2.21}
\end{align}
Since
$
\nlim\limits_{s\to s_0}H(s)=H(s_0),
$
one infers
\begin{equation}
\label{2.23}
\nlim_{s\to s_0}(H(s)-z)^{-1}=(H(s_0)-z)^{-1},
\end{equation}
uniformly with respect to $z\in \Gamma$. Thus,
combining \eqref{2.23} and \eqref{2.21}, one concludes
that $s\mapsto
\varphi(H(s))$ is differentiable with respect to $s\in S$ in
$\calB(\calH)$-topology. Along  with the existence of a 
continuous
$V''(s)$ on $S$, this yields the differentiability of the 
function in
\eqref{2.19}. Equation
\eqref{2.20} then follows by a straightforward computation
using \eqref{2.18}.
\end{proof}
%%%%%%%%%%%%%%%%%%%%%%%%%%%%%%%%%%%%%%%%%%%%%%%%%%%%%%%%%%%%%%%%%

%%%%%%%%%%%%%%%%%%%%%%%%%%%%%%%%%%%%%%%%%%%%%%%%%%%%%%%%%%%%%%%%%
\begin{lemma}\label{l3.4}
Under the assumptions of Lemma~\ref{l3.3} suppose in addition that
$V(s)$ is concave with respect to $s\in S$ in the sense that
\begin{equation}
\label{Vpp}
0\ge V''(s)\in \calB_1(\calH), \quad s\in S,
\end{equation}
and that the function $\varphi$ is nonnegative and
nonincreasing on $(a,b)$.  Then the function
\begin{equation}
s\longmapsto \tr\bigl(V'(s)\,\varphi(H(s))\bigr), \quad s\in S,
\end{equation}
is differentiable and  nonincreasing on $S$.
\end{lemma}
%%%%%%%%%%%%%%%%%%%%%%%%%%%%%%%%%%%%%%%%%%%%%%%%%%%%%%%%%%%%%%%%
\begin{proof}
Since $\varphi$ is nonnegative, $0 \le \varphi(H(s))\in
\calB(\calH)$, $s\in S$ and hence
\begin{equation}
\label{2.26}
\tr\bigl(V''(s)\,\varphi(H(s))\bigr)\le0, \quad s\in S,
\end{equation}
by \eqref{Vpp}.  Applying Theorem~\ref{t2.3} and
Remark~\ref{r2.4}\,(i) to the operator-valued
Herglotz functions
$M_1(z)=V'(s)(H(s)-z)^{-1}V'(s)$ and $M_2(z)=(H(s)-z)^{-1}$,
 one obtains the inequality
\begin{equation}
\label{old}
\frac{1}{2\pi i}\oint_{\Gamma}dz\,
\varphi(z)\,\tr\bigl[V'(s)(H(s)-z)^{-1}\bigr]^2\ge 0,
\quad s\in S.
\end{equation}
Combining \eqref{2.26} and \eqref{old} one infers
\begin{equation}
     \frac{d}{ds}\tr\bigl(V'(s)\,\varphi(H(s))\bigr)\le 0
\end{equation}
by Lemma~\ref{l3.3}, proving the assertion.
\end{proof}
%%%%%%%%%%%%%%%%%%%%%%%%%%%%%%%%%%%%%%%%%%%%%%%%%%%%%%%%%%%%%%%%%%%

In the general case of unbounded operators $H_0$ one can prove the
following result.
%%%%%%%%%%%%%%%%%%%%%%%%%%%%%%%%%%%%%%%%%%%%%%%%%%%%%%%%%%%%%%%%%%%
\begin{theorem}\label{t3.5}
Assume Hypothesis~\ref{h1.5}. Suppose in addition that $V'(s)$
is continuously differentiable in $\calB_1(\calH)$-norm
with respect to
$s\in (s_1, s_2)$ and that $V(s)$ is concave in the sense that
\begin{equation}
\label{4.2}
     0\ge V''(s) \in \calB_1(\calH), \quad s\in (s_1, s_2).
\end{equation}
Let $\varphi$ be a bounded nonnegative and  nonincreasing
real-analytic function on $\bbR$ admitting the analytic
continuation to a domain $\calD$ of the complex plane $\bbC$
containing the real axis $\bbR$.  Then the function
\begin{equation}
  s\longmapsto \tr\bigl(V'(s)\,\varphi(H(s))\bigr),
  \quad s\in (s_1,s_2),
\end{equation}
is nonincreasing with respect to $s\in (s_1,s_2)$.
\end{theorem}
%%%%%%%%%%%%%%%%%%%%%%%%%%%%%%%%%%%%%%%%%%%%%%%%%%%%%%%%%%%%%%%%
\begin{proof}
Introducing the sequence $\{P_n\}_{n\in\bbN}$ of
spectral projections of $H_0$,
\begin{equation}
   P_n=E_{H_{0}}((-n,n)), \quad n\in \bbN,
\end{equation}
one concludes that for fixed $s\in \bbR$ the bounded operators
given by
\begin{equation}
H^{(n)}(s)=P_nH_0P_n+V(s), \quad \dom(H^{(n)}(s))=\calH,
\end{equation}
converge to $H(s)$ in the strong resolvent sense and therefore,
\begin{equation}
\slim_{n\to \infty}\,\varphi(H^{(n)}(s))=\varphi(H(s)), \quad s\in
(s_1,s_2),
\end{equation}
by Theorem~VIII.20 in \cite{RS80}.

Since by hypothesis $\varphi$ is a bounded function, the family
of operators $\varphi(H^{(n)}(s))$ is uniformly bounded
with respect to
$n$ and hence by Theorem~1 in \cite{Gr73},
\begin{equation}
\label{3.6}
\lim_{n\to \infty}\,\tr\bigl(V'(s)
\,\varphi(H^{(n)}(s)\bigr)=
\tr\bigl(V'(s)\,\varphi(H(s))\bigr), \quad s\in (s_1,s_2).
\end{equation}

Given $n\in \bbN$ and $s_0\in (s_1,s_2)$, one can always find a
closed oriented Jordan contour $\Gamma$ in $\calD$, and a bounded
interval $(a,b)\subset\bbR$ such that the collection
$\{H^{(n)}(s_0),\calD,(a,b),\Gamma\}$ satisfies
Hypothesis~\ref{h3.1}. Since by hypothesis $\varphi$ is a 
nonnegative
nonincreasing function on $(a,b)$ and
\eqref{4.2} holds, one concludes by Lemma~\ref{l3.4} that the
function
\begin{equation}
\label{3.7}
   s\longmapsto \tr\bigl(V'(s)\,\varphi(H^{(n)}(s))\bigr)
\end{equation}
is nonincreasing in some neighborhood of $s_0$ and hence on
the whole interval $ (s_1,s_2)$ since $s_0\in (s_1,s_2)$ was
arbitrary.  Therefore, the function
\begin{equation}
    s\longmapsto \tr\bigl(V'(s)\,\varphi(H(s))\bigr),
     \quad s\in (s_1,s_2),
\end{equation}
is also nonincreasing on $ (s_1,s_2)$ as a pointwise limit
\eqref{3.6} of the nonincreasing functions in \eqref{3.7}.
\end{proof}
%%%%%%%%%%%%%%%%%%%%%%%%%%%%%%%%%%%%%%%%%%%%%%%%%%%%%%%%%%%%%%%

Now we are able to prove  Theorem~\ref{t1.7}, which is an 
extension 
of the monotonicity result, Theorem~\ref{t1.2}, of Birman and 
Solomyak \cite{BS75}.
 
\medskip

%%%%%%%%%%%%%%%%%%%%%%%%%%%%%%%%%%%%%%%%%%%%%%%%%%%%%%%%%%%%%%%
{\it Proof of Theorem~\ref{t1.7}.}\,
Given $\mu\in \bbR$, introducing the real-analytic function
\begin{equation}
\varphi_{\mu, \varepsilon}(\lambda)= \frac{1}{2}-\frac{1}
{\pi}\arctan \bigg
(\frac{\lambda-\mu}{\varepsilon}+\frac{1}{\sqrt{\varepsilon}}
\bigg), \quad \varepsilon >0,
\end{equation}
one concludes that
\begin{equation}
\label{3.11}
\lim_{\varepsilon \downarrow 0}\varphi_{\mu, 
\varepsilon}(\lambda)=
\chi_{(-\infty, \mu)}(\lambda),
\end{equation}
where $\chi_\Delta(\cdot)$ denotes the characteristic function
of the set $\Delta$.

Since
\begin{equation}
\label{3.12}
\sup_{\varepsilon>0}\|\varphi_{\mu, 
\varepsilon}\|_{L^\infty(\bbR)} <\infty,
\end{equation}
\eqref{3.11} implies the strong convergence
\begin{equation}
\label{3.13}
\slim_{\varepsilon \downarrow 0}\varphi_{\mu, 
\varepsilon}(H(s))=
E_{H(s)}((-\infty,\mu))
\end{equation}
by Theorem VIII.5 in \cite{RS80}.  Combining
\eqref{3.11}--\eqref{3.13} with Theorem~1 in \cite{Gr73},
one infers
\begin{equation}
\lim_{\varepsilon \downarrow 0}\,\tr\bigl(V'(s)
\,\varphi_{\mu, \varepsilon}(H(s))\bigr)=
\tr\bigl(V'(s)E_{H(s)}((-\infty,\mu))\bigr).
\end{equation}
By Theorem~\ref{t1.7}, the function
\begin{equation}
  s\longmapsto \tr\bigl(V'(s)\,\varphi_{\mu,
 \varepsilon}(H(s))\bigr),\quad s\in (s_1,s_2),
\end{equation}
is nonincreasing, proving the assertion, since the
pointwise limit of nonincreasing functions is nonincreasing.
\hfill $\square$
%%%%%%%%%%%%%%%%%%%%%%%%%%%%%%%%%%%%%%%%%%%%%%%%%%%%%%%%%%%%%

Next we prove Corollary~\ref{c1.8}.

\smallskip

%%%%%%%%%%%%%%%%%%%%%%%%%%%%%%%%%%%%%%%%%%%%%%%%%%%%%%%%%%%
{\it Proof of Corollary~\ref{c1.8}.}\,
By Theorem~\ref{t1.6}
\begin{equation}
\label{3.17}
\int_{-\infty}^\mu d\lambda\, \xi (\lambda , H_0, H(s))=
\int_{0}^{s} dt\,\tr\bigl(V'(t)E_{H(t)}((-\infty, \mu))\bigr),
 \quad s\in (s_1,s_2), \,\, \mu \in \bbR.
\end{equation}
By Theorem~\ref{t1.7} the integrand on the right-hand side
of \eqref{3.17}
is a nonincreasing function of $t$ and hence the
left-hand side of \eqref{3.17} is a concave function of 
$s$. Thus
\eqref{Concav} holds.

In order to prove \eqref{3.16a} one notes that
$\zeta(\mu,0+)=0$ and that a necessary and sufficient
condition for a measurable concave function $f(t)$  to be
subadditive on $(0,\infty)$ is that $f(0+)\ge 0$ (see, e.g.,
\cite[Theorem~7.2.5]{Hi48}).
\hfill $\square$
%%%%%%%%%%%%%%%%%%%%%%%%%%%%%%%%%%%%%%%%%%%%%%%%%%%%%%%%%%%%%

In the case of semibounded operators  the following statement 
might
be useful. We recall that $\calW_1(\bbR)$ denotes the 
function class
of $\varphi$ with  $\varphi'$  the Fourier transform of a 
finite
(complex) Borel measure (cf.~\eqref{5.2}).

%%%%%%%%%%%%%%%%%%%%%%%%%%%%%%%%%%%%%%%%%%%%%%%%%%%%%%%%%%%
\begin{theorem}\label{t3.6}
Let $H_0$ be a self-adjoint operator in $\mathcal H$, bounded 
from below, and assume the hypotheses of Theorem~\ref{t1.7}. 
Denote by $\Lambda$ the smallest semi-infinite interval 
containing the spectra of the family $H(s)$, $s\in (s_1, s_2)$,
\begin{equation}
\label{5.9}
\Lambda= \bigg[\inf_{s\in(s_1,s_2)}\spec(H(s)),\infty\bigg).
\end{equation}
Let $\varphi\in \calW_1(\bbR)\cap C^2(\bbR)$ be concave
on $\Lambda$ in the sense that
\begin{equation}
\label{5.10}
\varphi''\bigl|_{\Lambda} \bigr. \le 0,
\end{equation}
and
\begin{equation}
\label{5.11}
\varphi'(\lambda)=o(1) \text{ as }\lambda \to +\infty\,.
\end{equation}
Then the function
\begin{equation}
\label{vexx}
  s\longmapsto\tr\bigl[\varphi\bigl(H(s)\bigr)-
\varphi(H_0)\bigr],
    \quad  s\in (s_1,s_2),
\end{equation}
is concave in the sense that for any $s,t\in (s_1,s_2)$ and for
any $0\le \alpha \le 1$,
\begin{align}
&\lefteqn{\tr\bigl[\varphi(H(\alpha s+(1-\alpha)t))-
\varphi(H_0)\bigr]}
\no \\
&\geq\alpha \,\tr\bigl[\varphi(H(s))-\varphi(H_0)\bigr]+
(1-\alpha)\,\tr\bigl[\varphi(H(t))-\varphi(H_0)\bigr]\,.
\end{align}
In particular,
\begin{equation}
  \big[\varphi(H(\alpha s+(1-\alpha)t))- \alpha\,\varphi(H(s))-
  (1-\alpha)\,\varphi(H_{t})\big]\in \calB_1(\calH)
\end{equation}
and
\begin{equation}
   \tr \bigl[\varphi(H(\alpha s+(1-\alpha)t))- \alpha
   \,\varphi(H(s))- (1-\alpha)\,\varphi(H(s))\bigr]\ge 0.
\end{equation}
\end{theorem}
%%%%%%%%%%%%%%%%%%%%%%%%%%%%%%%%%%%%%%%%%%%%%%%%%%%%%%%%%%%%%%%%
\begin{proof}
First, one observes that by \eqref{ksi1} the integrated spectral
shift function $\zeta(\lambda, s)$ given by \eqref{ssf} is
uniformly bounded, that is,
\begin{equation}
|\zeta(\lambda, s)|\le \|V(s)\|_{\calB_1(\calH)},
\quad
\lambda\in \bbR, \,\,\, s\in (s_1,s_2).
\end{equation}
Moreover, since  $H_0$ is semibounded, one concludes that 
$\bbR\backslash\Lambda\ne\emptyset$ by definition \eqref{5.9}
of the set $\Lambda$ and hence for all $s\in (s_1,s_2)$,
\begin{equation}
\label{5.15}
\zeta(\lambda,s)=0, \quad \lambda\in \bbR\backslash \Lambda.
\end{equation}
Thus, using \eqref{5.11} and \eqref{5.15}, one infers
\begin{equation}
\label{5.16}
\lim_{\lambda\to \pm \infty}
\,\varphi'(\lambda)\,\zeta(\lambda,s)=0,\quad s\in (s_1, s_2)\,.
\end{equation}

Next, combining \eqref{5.15} and \eqref{5.16}, an integration
by parts in the trace formula \eqref{5.1} yields
\begin{equation}
\label{5.13}
\tr\bigl[\varphi(H(s))-\varphi(H_0)\bigr]=-\int_{\Lambda}
d\lambda\, \varphi''(\lambda)\,\zeta(\lambda,s).
\end{equation}
Given $\lambda \in \bbR$, the integrated spectral shift function
$\zeta(\lambda, s)$ is concave with respect to $s\in (s_1,s_2)$ by
Corollary~\ref{c1.8} and hence the left-hand side of \eqref{5.13} is
also a concave function of $s$  by \eqref{5.10} (as a weighted
mean of concave functions with a positive weight).
\end{proof}
%%%%%%%%%%%%%%%%%%%%%%%%%%%%%%%%%%%%%%%%%%%%%%%%%%%%%%%%%%%%%%%%%
\begin{remark} \lb{r3.7}
(i) If the measure $\nu$ in representation \eqref{5.2} is
absolutely continuous, then condition \eqref{5.11} holds
automatically by the Riemann--Lebesgue Lemma.

%%%%%%%%%%%%%%%%%%%%%%%%%%%%%%%%%%%%%%%%%%%%%%%%%%%%%%%%%%%%%%%%
\smallskip

\noindent (ii) If  $\varphi$ is  convex on $\Lambda$, that is,
\begin{equation}
\varphi''\bigl|_\Lambda \bigr. \ge 0,
\end{equation}
then the function  given by \eqref{vexx} is convex.
\end{remark}
%%%%%%%%%%%%%%%%%%%%%%%%%%%%%%%%%%%%%%%%%%%%%%%%%%%%%%%%%%%%%%%%
\begin{example} \lb{e3.8}
Under assumptions of Theorem~\ref{t3.6}, choosing
$\varphi\in\calW_1(\bbR)$ as
\begin{equation}
   \varphi(\lambda)=\exp\big({-\lambda t}\big), \quad \lambda\in
   \Lambda, \,\,t\ge0,
\end{equation}
one concludes that for any $t>0$,
\begin{equation}
   s\longrightarrow\tr
   \bigl[\exp\bigl(-tH(s)\bigr)-\exp\bigl(-tH_0\bigr)\bigr]
\end{equation}
is a convex function of $s\in (s_1,s_2)$.
\end{example}
%%%%%%%%%%%%%%%%%%%%%%%%%%%%%%%%%%%%%%%%%%%%%%%%%%%%%%

%%%%%%%%%%%%%%%%%%%%%%%%%%%%%%%%%%%%%%%%%%%%%%%%%%%%%%%%%%%
\vspace*{2mm}
\noindent {\bf Acknowledgments.} We are indebted to Vadim 
Kostrykin for kindly making available to us 
reference~\cite{Ko99} prior to its publication.
%%%%%%%%%%%%%%%%%%%%%%%%%%%%%%%%%%%%%%%%%%%%%%%%%%%%%%%%%%%

%%%%%%%%%%%%%%%%%%%%%%%%%%%%%%%%%%%%%%%%%%%%%%%%%%%%%%%

\end{document}